\documentclass[11pt]{amsart}

\usepackage[T1]{fontenc}
\usepackage{lmodern}
\usepackage[margin=0.85in]{geometry}
\usepackage{amsmath,amssymb,mathtools}
\usepackage{microtype}
\usepackage{graphicx}
\usepackage{tikz}
\usetikzlibrary{arrows.meta,calc}
\usepackage{float}
\usepackage{enumitem}
\usepackage{needspace}
\usepackage[colorlinks=true,linkcolor=blue!55!black,citecolor=blue!55!black,
  urlcolor=blue!55!black]{hyperref}

\newtheorem{theorem}{Theorem}
\newtheorem{proposition}[theorem]{Proposition}
\newtheorem{lemma}[theorem]{Lemma}
\newtheorem{corollary}[theorem]{Corollary}
\newtheorem{remark}[theorem]{Remark}
\newtheorem{question}[theorem]{Open question}
\theoremstyle{definition}
\newtheorem{definition}[theorem]{Definition}

\newcommand{\Binc}{B_{\mathrm{inc}}}
\newcommand{\ind}{\mathbf{1}}
\newcommand{\TS}{T_S}
\newcommand{\TLR}{T_{LR}}
\newcommand{\arxiv}[1]{\href{https://arxiv.org/abs/#1}{arXiv:#1}}
\newcommand{\leflip}{\le_{\mathrm{flip}}}
\newcommand{\nleflip}{\nleq_{\mathrm{flip}}}

\definecolor{agreen}{HTML}{2E7D32}
\definecolor{bred}{HTML}{C62828}
\definecolor{lightgreen}{HTML}{DFF2DF}
\definecolor{lightred}{HTML}{F6DADA}
\definecolor{alight}{HTML}{E5F3E6}
\definecolor{blight}{HTML}{F8E4E2}
\definecolor{gridgray}{HTML}{D8D8D8}
\definecolor{muted}{HTML}{666666}
\definecolor{splitgray}{HTML}{888888}
\newcommand{\diagramfont}{\normalfont\small}
\newcommand{\diagramdetailfont}{\normalfont\scriptsize}
\newcommand{\pathpanel}[6]{%
  \begin{scope}[shift={(#1,#2)}]
    \node[font=\diagramfont,text=#3] at (2,4.95) {#4};
    \node[font=\diagramdetailfont,text=muted] at (2,4.52) {#5};
    \fill[#3!10] (-.18,-.18) rectangle (4.18,4.18);
    \draw[gridgray,line width=.26475pt] (0,0) grid[step=1] (4,4);
    \draw[splitgray,dashed,line width=.45386pt] (0,4)--(4,0);
    #6
  \end{scope}%
}
\newcommand{\stepnum}[4]{\node[font=\diagramdetailfont,text=#4] at (#1,#2) {#3};}
\newcommand{\bit}[4]{\node[font=\diagramdetailfont,text=#4,anchor=base] at (#1,#2) {#3};}
\colorlet{included}{green!55!black}
\colorlet{excluded}{red!80!black}

\hypersetup{
  pdftitle={Poset dimensions of higher Bruhat orders},
  pdfsubject={Order dimension for inclusion and flip orders on higher Bruhat inversion sets}
}
\title[Dimension bounds for higher Bruhat orders]{Dimension bounds for higher Bruhat orders}
\date{}
\author{Daria Poliakova}
\address{Universit\"at Hamburg, Hamburg, Germany}
\email{polydarya@gmail.com}

\begin{document}
\raggedbottom

\begin{abstract}
We study the Dushnik--Miller dimension of higher Bruhat orders. For the inclusion order
$\Binc(n,k)$, with $k$ fixed, we prove that the order of dimension growth is $n^k$, by exhibiting standard examples for the lower bound and explicit systems of coordinates for the upper bound. A recent
Tur\'an-system construction of Gao, Kuang, Pikhurko and Wang gives an upper
coefficient within a factor smaller than $4$ of the lower
coefficient. For the genuine higher Bruhat order $B(n,k)$, we construct for every fixed
$k\ge3$ a standard example of order $n^{2k+1}$. We leave it open whether a polynomial upper bound on this dimension exists.
\end{abstract}

\maketitle

\section*{Introduction}

Dushnik--Miller dimension is a well-known measure of poset complexity: for a poset \(P\), \(\dim P\) is the smallest \(d\) such that \(P\) admits an order-embedding into a product of \(d\) chains \cite{DushnikMiller}. Poset dimensions are known for many interesting posets. For the weak order on the symmetric group \(S_n\), one has \(\dim S_n=n-1\) \cite{Flath1993,reading}. The same equality holds in types \(B\) and \(D\) \cite{reading}, and has recently been shown by me \cite{PoliakovaRank} and independently by Segovia to fail in some of the exceptional types. 

Higher Bruhat orders \cite{ManinSchechtman,Ziegler} are another famous generalization of the weak order on the symmetric group. Instead of considering orderings of elements of \([n]\), one considers appropriately defined equivalence classes of admissible orderings of \(k\)-element subsets of \([n]\). The inversion sets of such orderings are precisely the consistent sets of \((k+1)\)-element subsets of \([n]\) \cite{Ziegler}. Unlike the \(k=1\) case, for general \(k\) the order of such inversion sets by ordinary inclusion is stronger than the flip order, where one is only allowed to add one inversion at a time: flip comparability implies inclusion, but the converse may fail. 

Somewhat surprisingly, to our knowledge the dimension theory of higher Bruhat orders beyond the \(k=1\) case has never been studied. The current paper aims to close this gap in the literature. Computer experiments indicate that the precise values of \(\dim\Binc(n,k)\) and \(\dim B(n,k)\) are really difficult to describe, so we lower our expectations and ask an asymptotic question: for \(k\) fixed, how do \(\dim\Binc(n,k)\) and \(\dim B(n,k)\) grow with \(n\)? Our main results are
\[
\dim\Binc(n,k)=\Theta_k(n^k)\qquad (k\geq2),
\]
and
\[
\dim B(n,k)=\Omega_k(n^{2k+1})\qquad (k\geq3).
\]

For the lower bound on the dimension of the inclusion order, we construct a Boolean interval, using ideas of Chau \cite{Chau2024}. For the upper bound, we introduce a system of coordinates generalizing the usual weighted Lehmer codes \cite{reading,poliakova}. A very recent development in Tur\'an systems due to Gao, Kuang, Pikhurko and Wang \cite{GKPW} allows us to choose a smaller coordinate subsystem sufficient to detect all failures of inclusion, making the ratio of the leading coefficients of \(n^k\) in the lower and upper bounds independent of \(k\), and in fact smaller than \(4\).

We then demonstrate that the flip order on \(B(n,k)\) is measurably more complex: we construct a standard example of size \(\Omega_k(n^{2k+1})\), generalizing Ziegler's example of inclusion failure in \(B(8,3)\) \cite[Theorem~4.5]{Ziegler}. In this case, we make no claim of optimality for our lower bound: we did not succeed in obtaining any polynomial upper bound on \(\dim B(n,k)\). We leave the question of whether such a bound exists open. 

A different recent manifestation of the complexity of the flip order is
that computing flip distance in $B(n,k)$ is NP-complete for every fixed
$k\geq2$ \cite{Marc2026}.

{\bf Acknowledgements.} I am grateful to Vadim Schechtman and Gleb Koshevoy for discussing higher Bruhat orders with me, and to Andrii Bondarenko for teaching me how to ask asymptotic questions. I was supported by the Deutsche Forschungsgemeinschaft (DFG, German Research Foundation) -- SFB-Gesch\"aftszeichen 1624 -- Projektnummer 506632645.

\section{Preliminaries}
Fix integers $k\ge2$ and $n\ge k+1$, and write $[n]=\{1,\ldots,n\}$.

\begin{definition}\label{def:consistent}
For $M=\{m_1<\cdots<m_{k+2}\}\subseteq[n]$, its \emph{packet} is the ordered list
\[
P(M)=\bigl(M\setminus\{m_{k+2}\},M\setminus\{m_{k+1}\},\ldots,M\setminus\{m_1\}\bigr).
\]
A family $X\subseteq\binom{[n]}{k+1}$ is \emph{consistent} if its members in each
packet form a beginning or an ending segment; the empty segment and the whole
packet are allowed. 
\end{definition}

On the set of consistent families, we consider two orders. The poset $\Binc(n,k)$ is ordered by inclusion. The order in the poset
$B(n,k)$ is generated by consistent
single-set additions; we write $U\leflip V$ for this order. Thus
$U\leflip V$ implies $U\subseteq V$, but the converse may fail.

We also use an alternative characterization of consistent families. For a consistent family $W\subseteq\binom{[n]}{k+2}$, let $\preceq_W$ be the transitive closure of the packet orders, read in the usual direction for $H\notin W$ and in reverse for $H\in W$. Then $X\subseteq\binom{[n]}{k+1}$ is consistent if and only if it is an order ideal of $\preceq_W$ for some consistent $W$ \cite[Theorem~4.1]{Ziegler}.

For a poset $P$, its Dushnik--Miller dimension $\dim P$ is the least number of chains whose product contains an isomorphic
copy of $P$ \cite[Chapter~1]{Trotter}. 

\section{The inclusion order}

\subsection{Lower bounds via Boolean intervals}

For a set \(S\), the Boolean poset of its subsets \(2^S\) has dimension \(|S|\) \cite[Chapter~1]{Trotter}. Therefore any Boolean interval in \(\Binc(n,k)\) gives a lower bound on dimension. We start from Boolean intervals of Chau \cite{Chau2024}, and then explain a certain modular enlargement of those.

For a finite set $A$ of integers, put $\sigma(A)=\sum_{a\in A}a$, and consider level sets $S_c=\{A\in\binom{[n]}{k+1}:\sigma(A)=c\}$ and sublevel sets $X_c=\{A\in\binom{[n]}{k+1}:\sigma(A)<c\}$. Sums strictly increase when moving along one packet; thus, for \(c\) fixed, $X_c$ is a possibly empty beginning segment in each packet and $S_c$ intersects this packet at most once, right after this beginning segment. It follows that both $X_c$ and $X_c\cup E$ for any $E\subseteq S_c$ are consistent, so $[X_c,X_c\cup S_c]\simeq 2^{S_c}$ is a Boolean interval of dimension $|S_c|$. These are the Boolean intervals implicit in Chau's construction \cite[Lemma~4.8 and (6.2)--(6.6)]{Chau2024}.

We then obtain a naive lower bound on dimension by the pigeonhole principle: possible sums $c$ range from $1+\cdots+(k+1)=\binom{k+2}{2}$ to $(n-k)+\cdots+n=(k+1)n-\binom{k+1}{2}$, so there exists $S_c$ with
$|S_c|\ge \binom{n}{k+1}/((k+1)(n-k-1)+1)$. Expanding in $n$ for fixed $k$, we get
\[
\dim\Binc(n,k)\geq \frac{1}{(k+1)(k+1)!}n^k+O_k(n^{k-1}).
\]
By a finer estimate on the size of the appropriate level set \cite[(6.3)--(6.6)]{Chau2024}, Chau improves the naive constant $1/((k+1)(k+1)!)$ to
$A(k,\lfloor (k-1)/2\rfloor)/(k!(k+1)!)$, where $A(k,j)$ is the Eulerian number counting permutations of $[k]$ with exactly $j$ descents. We now give a stronger improvement by passing to modular level sets. Under the identification of $(k+1)$-subsets of $[n]$ with binary words of weight $k+1$, our modular level sets are precisely the Graham--Sloane constant-weight codes of minimum distance $4$ \cite[Theorem~1]{GrahamSloane1980}.

\begin{proposition}\label{prop:lower}
For $k\ge 2$ and $n\ge k+1$,
\[
\dim\Binc(n,k)\ge\left\lceil\frac1n\binom{n}{k+1}\right\rceil.
\]
\end{proposition}

\begin{proof}
Let $c$ be a residue modulo $n$. Set the modular level set
$\bar{S}_c=\{A\in\binom{[n]}{k+1}:\sigma(A)\equiv c\pmod n\}$.
Choose an integer representative of $c$ and put
\[
\bar{X}_c=\left\{A\in\binom{[n]}{k+1}:
 c+jn<\sigma(A)<c+(j+1)n\text{ for some even }j\in\mathbb Z\right\}.
\]

\begin{figure}[H]
\centering
\begin{tikzpicture}[x=1.8cm,y=1cm,font=\diagramfont]
  \foreach \j in {-3,-2,-1,0,1,2} {
    \pgfmathtruncatemacro{\parity}{mod(\j+4,2)}
    \ifnum\parity=0
      \def\intervalcolour{included}
    \else
      \def\intervalcolour{excluded}
    \fi
    \fill[\intervalcolour!13] (\j,-0.07) rectangle ({\j+1},0.07);
    \draw[\intervalcolour,line width=1.9pt]
      (\j,0)--({\j+1},0);
    \draw[<->,>=Stealth,line width=0.35pt]
      (\j,0.40)--({\j+1},0.40)
      node[midway,above=1pt] {$n$};
    \node at ({\j+0.5},0.98) {$j=\j$};
  }
  \foreach \x/\lab in {-3/{c-3n},-2/{c-2n},-1/{c-n},0/{c},
                       1/{c+n},2/{c+2n},3/{c+3n}} {
    \filldraw[fill=white,draw=black,line width=0.65pt]
      (\x,0) circle[radius=2.4pt];
    \node[below=7pt,font=\diagramfont] at (\x,0) {$\lab$};
  }
  \node at (-3.36,0) {$\cdots$};
  \node at (3.33,0) {$\cdots$};
  \draw[->,>=Stealth,line width=0.4pt] (3.60,0)--(3.97,0);
  \node[below=7pt,font=\diagramfont] at (3.77,0) {$\sigma(A)$};
\end{tikzpicture}
\caption{Green intervals give $\bar{X}_c$; red intervals stay excluded.
White points correspond to elements of $\bar{S}_c$.}
\label{fig:intervals}
\end{figure}
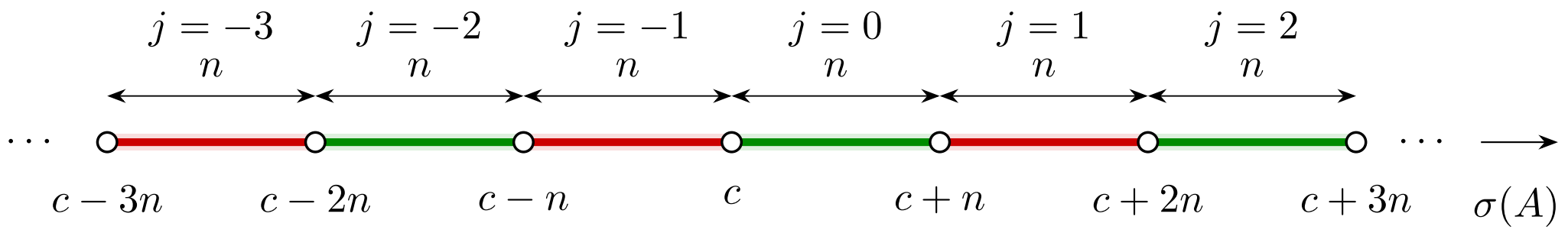

Along each packet, the sums strictly increase and span an interval of length
$\max M-\min M<n$. As Figure~\ref{fig:intervals} shows, sampling colors at these
sums gives at most one color change along a packet, no matter how the possible white point is colored.
Thus $\bar{X}_c\cup E$ is consistent for every $E\subseteq\bar{S}_c$, so
$[\bar{X}_c,\bar{X}_c\cup\bar{S}_c]\cong 2^{\bar{S}_c}$ is a Boolean interval.

The $n$ modular level sets $\bar{S}_c$ partition $\binom{[n]}{k+1}$, so by the pigeonhole principle some $c$ satisfies
$|\bar{S}_c|\ge\lceil n^{-1}\binom{n}{k+1}\rceil$. The corresponding Boolean interval has this dimension, finishing the proof.
\end{proof}

Expanding in \(n\) for fixed \(k\), we obtain
\[
\dim\Binc(n,k)\geq \frac{1}{(k+1)!}n^k-O_k(n^{k-1}).
\]

Thus for expansions in \(n\), the naive leading constant was $1/((k+1)(k+1)!)$, Chau's larger leading constant was
$A(k,\lfloor (k-1)/2\rfloor)/(k!(k+1)!)$, and our yet larger leading constant is
$1/(k+1)!$.

\begin{remark}
The same construction also slightly improves the counting lower bound of Chau. Indeed, since $\Binc(n,k)$ and $B(n,k)$ have the same ground set, the Boolean interval above gives at least $2^{|\bar S_c|}$ distinct consistent families. Hence
\[
\log_2 |B(n,k)|
\ge \frac1n\binom{n}{k+1}
= \frac{1}{(k+1)!}n^k-O_k(n^{k-1}).
\]
This improves the leading constant $A(k,\lfloor(k-1)/2\rfloor)/(k!(k+1)!)$ in Chau's lower bound \cite[Theorem~6.1]{Chau2024} to $1/(k+1)!$.
\end{remark}

\subsection{Long left and right coordinates}

The very naive upper bound on order dimension uses the indicator vector
$(\ind_{A\in X})_{A\in\binom{[n]}{k+1}}$.
It gives an order embedding into $\{0,1\}^{\binom{n}{k+1}}$, and hence
$\dim\Binc(n,k)\le\binom{n}{k+1}$. This is of order $n^{k+1}$.

\begin{definition}\label{def:coordinates}
For each $Q\in\binom{[n]}k$, define two binary words by
\[
\rho_Q(X)=
 \bigl(\ind_{Q\cup\{x\}\in X}\bigr)_{
    \substack{x\in[n]\setminus Q,\ x>\min Q\\x\text{ increasing}}},
\qquad
\lambda_Q(X)=
 \bigl(\ind_{Q\cup\{x\}\in X}\bigr)_{
    \substack{x\in[n]\setminus Q,\ x<\max Q\\x\text{ decreasing}}}.
\]
These are the \emph{long right} and \emph{long left coordinates}, respectively.
For different $X$, compare the corresponding words lexicographically.
\end{definition}

In other words, the right coordinate reads possible added elements from left
to right, starting just after $\min Q$; the left coordinate reads them from
right to left, starting just before $\max Q$. Elements already in $Q$ are
skipped. Each coordinate is order-preserving under inclusion.

The prototype of these long coordinates is Reading's coordinates for the weak order \cite{reading}, corresponding to the $k=1$ case. In that case, $Q = \{q\}$ is just one number, and $\lambda_Q(X)$ is a binary word encoding all inversions with $q$; see also \cite{poliakova} for the binary word description of Reading's coordinates.

\begin{lemma}\label{lem:safe}
Let $X,Y\in\Binc(n,k)$ with $X\nsubseteq Y$, and choose
$A=\{a_1<\cdots<a_{k+1}\}\in X\setminus Y$ with its span $a_{k+1}-a_1$
as small as possible. For $Q=A\setminus\{a_i\}$,
\begin{enumerate}[label=(\roman*),leftmargin=2.2em,itemsep=1pt,topsep=3pt]
\item if $i\ge2$, then $\rho_Q(X)>\rho_Q(Y)$;
\item if $i\le k$, then $\lambda_Q(X)>\lambda_Q(Y)$.
\end{enumerate}
\end{lemma}

\begin{proof}
Suppose $i\ge2$ and $\rho_Q(X)\le\rho_Q(Y)$.
At $x=a_i$, the two words have entries $1$ and $0$.
An earlier differing entry must therefore give $a_1<b<a_i$ with
$Q\cup\{b\}\in Y\setminus X$. Put $M=A\cup\{b\}$.
In $P(M)$, the member $Q\cup\{b\}=M\setminus\{a_i\}$ goes before
$A=M\setminus\{b\}$. Since $X$ excludes the first of these and includes
the second, its members form an ending segment. The members of $Y$ form a
beginning segment. The last packet member $M\setminus\{a_1\}$ therefore
belongs to $X\setminus Y$. Its minimum is larger than $a_1$ and its maximum
is $a_{k+1}$, so its span is smaller, contradicting the choice of $A$.

For the left coordinate, suppose $i\le k$ and
$\lambda_Q(X)\le\lambda_Q(Y)$. An earlier differing entry in the decreasing
word gives $a_i<b<a_{k+1}$ with $Q\cup\{b\}\in Y\setminus X$.
Put $M=A\cup\{b\}$. Now $A=M\setminus\{b\}$ precedes
$Q\cup\{b\}=M\setminus\{a_i\}$ in $P(M)$.
The members of $X$ form a beginning segment and those of $Y$ an ending
segment. The first packet member $M\setminus\{a_{k+1}\}$ belongs to
$X\setminus Y$ and has smaller span than $A$, again a contradiction.
\end{proof}

\begin{corollary}\label{cor:allright}
Taking all right coordinates gives $\dim\Binc(n,k)\le\binom nk$.
\end{corollary}

\begin{proof}
If $X\nsubseteq Y$, take $A$ as in Lemma~\ref{lem:safe} and delete its
maximum. The resulting right coordinate has a larger value on $X$ than on
$Y$. Thus the coordinates give an order embedding into a product of
$\binom nk$ chains.
\end{proof}

Together, Proposition~\ref{prop:lower} and Corollary~\ref{cor:allright} give
$\dim\Binc(n,k)=\Theta_k(n^k)$. The leading coefficients are $1/(k+1)!$ and
$1/k!$, with ratio $k+1$. This makes us ask: can we choose fewer coordinates so that this ratio
is bounded independently of $k$? Below we explain the (positive) answer.

\subsection{Colored Tur\'an systems}

We use colouring templates as in \cite[Definition~3.11 and Section~3.3]{FOU},
allowing empty colour sets as in \cite[Definition~1.1]{FMR}.
We introduce colored Tur\'an systems.

\begin{definition}\label{def:colored}
Let $\Gamma$ be a finite set of colours. A \emph{colouring template} on
$\binom{[n]}k$ assigns to each $k$-set $Q$ a set $C(Q)\subseteq\Gamma$,
possibly empty (alternatively: for each colour independently, we choose a family of $k$-sets of this colour). The size of a colouring template is $|C|=\sum_{Q\in\binom{[n]}k}|C(Q)|$.
Fix the \emph{colouring rule} $S=(S_1,\ldots,S_{k+1})$, where each nonempty set $S_i\subseteq\Gamma$
specifies the colours allowed to cover a $(k+1)$-set through the facet
obtained by deleting its $i$-th element.
A colouring template $C$ is an \emph{$S$-colored Tur\'an system} if every
$A=\{a_1<\cdots<a_{k+1}\}\subseteq[n]$ satisfies
\[
C(A\setminus\{a_i\})\cap S_i\ne\varnothing
\quad\text{for at least one }i\in\{1,\ldots,k+1\}.
\]
\end{definition}

Write $\TS(n,k+1,k)$ for the minimum size of an $S$-colored Tur\'an system
on $\binom{[n]}k$.

For our coordinates, $\Gamma=\{L,R\}$ and the colouring rule is given by $S_1=\{L\}$,
$S_{k+1}=\{R\}$, and $S_i=\{L,R\}$ for $2\le i\le k$.
We call these \emph{$LR$-colored Tur\'an systems} and write
$\TLR(n,k+1,k)$ for their minimum size.
Thus deleting the minimum allows $L$, deleting the maximum allows $R$,
and deleting an interior element allows either colour.

\Needspace{7\baselineskip}
\begin{proposition}\label{prop:embedding}
Every $LR$-colored Tur\'an system $C$ gives an order embedding of
$\Binc(n,k)$ into a product of $|C|$ chains. Hence
$\dim\Binc(n,k)\le\TLR(n,k+1,k)$.
\end{proposition}

\begin{proof}
Keep $\lambda_Q$ whenever $L\in C(Q)$ and $\rho_Q$ whenever $R\in C(Q)$.
These are $|C|$ order-preserving coordinates. If $X\nsubseteq Y$, choose
$A$ as in Lemma~\ref{lem:safe}. The $LR$-colored Tur\'an system $C$ supplies a facet of $A$ carrying an
allowed colour. By Lemma~\ref{lem:safe}, the corresponding coordinate has a larger
value on $X$ than on $Y$.
\end{proof}

Recall that $T(n,k+1,k)$ is the minimum size of a family of $k$-sets that
contains a facet of every $(k+1)$-set. Such a family is a \emph{Tur\'an
$(n,k+1,k)$-system}. With one colour, the systems in
Definition~\ref{def:colored} are exactly ordinary Tur\'an systems.

\begin{lemma}\label{lem:factor}
Let $h$ be the least size of a set of colours that intersects every $S_i$.
Then
\[
T(n,k+1,k)\le\TS(n,k+1,k)\le h\,T(n,k+1,k).
\]
In the $LR$ case, $h=2$.
\end{lemma}

\begin{proof}
Forgetting the colours of an $S$-colored Tur\'an system $C$ gives the
ordinary Tur\'an system $\{Q:C(Q)\ne\varnothing\}$, whose size is at most
$|C|$. Conversely, take an ordinary Tur\'an system of minimum size and give
each of its members the same $h$ colours meeting every $S_i$.
Whichever facet covers a given $(k+1)$-set, one of these colours is allowed
there.
In the $LR$ case, the first and last positions require $L$ and $R$
respectively, and these two colours suffice.
\end{proof}

By Proposition~\ref{prop:embedding} and Lemma~\ref{lem:factor}, an ordinary
Tur\'an bound of the form
\begin{equation}\label{eq:ordinary}
T(n,k+1,k)\le\frac{\gamma}{k+1}\binom nk,
\qquad\gamma\text{ independent of }n\text{ and }k,
\end{equation}
would give an upper leading coefficient $2\gamma/(k+1)!$ for dimension.
Its ratio to the lower leading coefficient would be at most $2\gamma$.

Pikhurko's 2025 work gave the first bound of the form~\eqref{eq:ordinary}
with an absolute constant.

\begin{theorem}[Pikhurko, {\cite[Theorem~1.1]{Pikhurko}}]\label{thm:pikhurko}
For all integers $n>k\ge1$, $T(n,k+1,k)\le\frac{6.239}{k+1}\binom nk$.
For sufficiently large $k$, the constant $6.239$ can be replaced by $4.911$,
uniformly for all $n>k$.
\end{theorem}

Thus Theorem~\ref{thm:pikhurko} already gives a ratio of leading coefficients
at most $12.478$, and at most $9.822$ for sufficiently large $k$.

\subsection{The Gao--Kuang--Pikhurko--Wang construction}

In very recent work, Gao, Kuang, Pikhurko and Wang
\cite[Theorem~1.1]{GKPW} reduced the factor in the ordinary Tur\'an bound to
$4(k+1)/(k+2)<4$.
We are lucky that their construction \cite[Section~3]{GKPW} can be read as
first selecting facets with colours $L$ and $R$, then forgetting the colours
to obtain an ordinary Tur\'an system. Retaining them gives an
$LR$-colored Tur\'an system directly, removing the factor of two in
Lemma~\ref{lem:factor}. We record their construction below, for the reader's convenience.

For a finite ordered set $S$ and $0\le j\le |S|$, write $L_j(S)$ and $R_j(S)$
for its first and last $j$ elements, respectively.

\begin{definition}\label{def:random}
Assume $k\ge3$ and put $m=\lfloor(k+1)/2\rfloor$.
Choose a total order $\prec$ on $\binom{[n]}{m-1}$.
Define a colouring template $C$ by assigning colours to each
$Q\in\binom{[n]}k$ as follows:
\begin{align}
&\text{assign }L\quad\text{ if }\quad
L_{m-1}(Q)\prec Z\quad\text{for every }Z\in\binom{R_m(Q)}{m-1},
\label{eq:left}\\
&\text{assign }R\quad\text{ if }\quad
R_{m-1}(Q)\prec Z\quad\text{for every }Z\in\binom{L_m(Q)}{m-1}.
\label{eq:right}
\end{align}
\end{definition}

In terms of the long coordinates $\lambda_Q$ and $\rho_Q$, this means choosing
$\lambda_Q$ when the set of the first $m-1$ elements of $Q$ precedes every
facet of its last $m$ elements, and choosing $\rho_Q$ when the set of its last
$m-1$ elements precedes every facet of its first $m$ elements.
Our parameters $k,m$ correspond to $r,k$ in \cite[Definition~3.1]{GKPW}.

\Needspace{7\baselineskip}
\begin{lemma}\label{lem:cover}
The template in Definition~\ref{def:random} is an $LR$-colored Tur\'an system
for every choice of $\prec$.
\end{lemma}

\begin{proof}
The argument in \cite[proof of Lemma~3.2]{GKPW} supplies in every
$(k+1)$-set $A$ a facet satisfying~\eqref{eq:left} by deleting one of the
first $m$ elements, or a facet satisfying~\eqref{eq:right} by deleting one of
the last $m$ elements. These deletions avoid $\max A$ and $\min A$,
respectively, so the assigned colour is allowed.
\end{proof}

\begin{lemma}\label{lem:expectation}
If $\prec$ is chosen uniformly among all total orders in
Definition~\ref{def:random}, then
\[
\mathbb E|C|=\frac{2}{m+1}\binom nk.
\]
\end{lemma}

\begin{proof}
The computation in \cite[proof of Lemma~3.3]{GKPW} gives
$\Pr(L\in C(Q))=\Pr(R\in C(Q))=1/(m+1)$ for every
$Q\in\binom{[n]}k$. Summing these probabilities over $Q$ gives the formula.
\end{proof}

\begin{theorem}\label{thm:main}
For $k\ge2$ and $n\ge k+1$,
\[
\left\lceil\frac1n\binom{n}{k+1}\right\rceil
\le\dim\Binc(n,k)\le\TLR(n,k+1,k)\le c_k\binom nk,
\]
where
\[
c_k=\frac{2}{1+\lfloor(k+1)/2\rfloor}\le\frac{4}{k+2}.
\]
For fixed $k$, these lower and upper bounds have leading coefficients
$1/(k+1)!$ and $c_k/k!$; their ratio is $(k+1)c_k<4$.
\end{theorem}

\begin{proof}
For $k\ge3$, Lemmas~\ref{lem:cover} and~\ref{lem:expectation} give an
$LR$-colored Tur\'an system of size at most
$\frac{2}{m+1}\binom nk=c_k\binom nk$.
For $k=2$, assigning $R$ to every $2$-set gives an $LR$-colored Tur\'an
system of size $\binom n2$, and $c_2=1$.
The dimension bounds follow from Propositions~\ref{prop:lower}
and~\ref{prop:embedding}. Finally,
$(k+1)c_k\le4(k+1)/(k+2)<4$ gives the claimed ratio.
\end{proof}

As a side comment, we believe that the true order of growth is closer to the lower bound than to the upper bound.

\section{The flip order}

We now turn to the genuine higher Bruhat order. Throughout this section
$k\ge3$. A \emph{flip path} from $U$ to $V$ is a sequence of consistent
single-set additions. We call $(U,V)$ a \emph{gap} if $U,V$ are consistent,
$U\subseteq V$, and $U\nleflip V$. Its two members are incomparable in the
flip order, since upward flips only add sets. 

For equally sized sets $F=\{f_1<\cdots<f_s\}$ and
$G=\{g_1<\cdots<g_s\}$, write $F\preceq G$ when $f_i\le g_i$ for every
$i$. A subset of a coordinatewise poset is a down-set or up-set in the usual
sense. We will construct large standard examples to prove the following lower bounds.

\begin{theorem}\label{thm:flip-main}
Let $k\ge3$ and $n\ge2k+2$, and put
$\ell=\lfloor n/2\rfloor$ and $m=\lceil n/2\rceil$. Then
\[
\dim B(n,k)\ge
\left\lceil
\frac{\binom{\ell}{k+1}\binom{m}{k+1}}
{(k+1)(n-2k-2)+1}
\right\rceil.
\]
Consequently, for fixed $k$,
\[
\dim B(n,k)\ge
\frac{n^{2k+1}}{2^{2k+2}(k+1)((k+1)!)^2}-O_k(n^{2k}).
\]
\end{theorem}

\subsection{The local gap}
In this subsection, we generalize Ziegler's example of a gap in $B(8,3)$ to any $B(2(k+1),k)$.
\begin{definition}\label{def:local-gap}
Take
\[
A=\{a_1<\cdots<a_{k+1}\}<B=\{b_1<\cdots<b_{k+1}\},
\]
where every point of $A$ is smaller than every point of $B$. A $(k+1)$-subset of
$A\cup B$ is \emph{mixed} if it meets both $A$ and $B$. A mixed set is
\emph{$A$-leaning} if it avoids $b_k,b_{k+1}$ and either contains $a_1,a_2$
or has only one point in $B$. It is \emph{$B$-leaning} if it avoids
$a_1,a_2$ and either contains $b_k,b_{k+1}$ or has only one point in $A$.
Put
\[
X_{AB}=\{A,B\},\qquad
Y_{AB}=\left\{F\in\binom{A\cup B}{k+1}:F\text{ is neither $A$- nor $B$-leaning}\right\}.
\]
For $k=3$, this is exactly Ziegler's gap in $B(8,3)$
\cite[Theorem~4.5]{Ziegler}.
\end{definition}

\begin{figure}[H]
\centering
\begin{tikzpicture}[scale=.756437,x=1.15cm,y=1.15cm,font=\diagramfont,line cap=round,line join=round,line width=.30257pt]
  \node[rotate=90,font=\diagramfont,text=agreen] at (-.95,8.70) {$A$-leaning};
  \node[rotate=90,font=\diagramfont,text=bred] at (-.95,2.20) {$B$-leaning};

  \pathpanel{0}{6.7}{agreen}{\textcolor{agreen}{12}\textcolor{bred}{56}}{$1,2$ chosen}{
    \draw[agreen,line width=1.66416pt] (0,0)--(0,1)--(0,2)--(1,2)--(2,2)--(2,3)--(2,4)--(3,4)--(4,4);
    \fill (0,0) circle(1.4pt); \fill (4,4) circle(1.4pt);
    \fill[white] (2,2) circle(2.1pt); \draw (2,2) circle(2.1pt);
    \stepnum{.18}{.5}{1}{agreen}\stepnum{.18}{1.5}{2}{agreen}\stepnum{.5}{2.18}{3}{agreen}\stepnum{1.5}{2.18}{4}{agreen}
    \stepnum{2.18}{2.5}{5}{bred}\stepnum{2.18}{3.5}{6}{bred}\stepnum{2.5}{3.82}{7}{bred}\stepnum{3.5}{3.82}{8}{bred}
    \bit{1.0}{-.55}{1}{agreen}\bit{1.3}{-.55}{1}{agreen}\bit{1.6}{-.55}{0}{agreen}\bit{1.9}{-.55}{0}{agreen}
    \bit{2.35}{-.55}{1}{bred}\bit{2.65}{-.55}{1}{bred}\bit{2.95}{-.55}{0}{bred}\bit{3.25}{-.55}{0}{bred}
  }

  \pathpanel{5.05}{6.7}{agreen}{\textcolor{agreen}{234}\textcolor{bred}{5}}{one $B$-point}{
    \draw[agreen,line width=1.66416pt] (0,0)--(1,0)--(1,1)--(1,2)--(1,3)--(1,4)--(2,4)--(3,4)--(4,4);
    \fill (0,0) circle(1.4pt); \fill (4,4) circle(1.4pt);
    \fill[white] (1,3) circle(2.1pt); \draw (1,3) circle(2.1pt);
    \stepnum{.5}{.18}{1}{agreen}\stepnum{1.18}{.5}{2}{agreen}\stepnum{1.18}{1.5}{3}{agreen}\stepnum{1.18}{2.5}{4}{agreen}
    \stepnum{1.18}{3.5}{5}{bred}\stepnum{1.5}{3.82}{6}{bred}\stepnum{2.5}{3.82}{7}{bred}\stepnum{3.5}{3.82}{8}{bred}
    \bit{1.0}{-.55}{0}{agreen}\bit{1.3}{-.55}{1}{agreen}\bit{1.6}{-.55}{1}{agreen}\bit{1.9}{-.55}{1}{agreen}
    \bit{2.35}{-.55}{1}{bred}\bit{2.65}{-.55}{0}{bred}\bit{2.95}{-.55}{0}{bred}\bit{3.25}{-.55}{0}{bred}
  }

  \pathpanel{10.1}{6.7}{agreen}{\textcolor{agreen}{123}\textcolor{bred}{5}}{both mechanisms}{
    \draw[agreen,line width=1.66416pt] (0,0)--(0,1)--(0,2)--(0,3)--(1,3)--(1,4)--(2,4)--(3,4)--(4,4);
    \fill (0,0) circle(1.4pt); \fill (4,4) circle(1.4pt);
    \fill[white] (1,3) circle(2.1pt); \draw (1,3) circle(2.1pt);
    \stepnum{.18}{.5}{1}{agreen}\stepnum{.18}{1.5}{2}{agreen}\stepnum{.18}{2.5}{3}{agreen}\stepnum{.5}{3.18}{4}{agreen}
    \stepnum{1.18}{3.5}{5}{bred}\stepnum{1.5}{3.82}{6}{bred}\stepnum{2.5}{3.82}{7}{bred}\stepnum{3.5}{3.82}{8}{bred}
    \bit{1.0}{-.55}{1}{agreen}\bit{1.3}{-.55}{1}{agreen}\bit{1.6}{-.55}{1}{agreen}\bit{1.9}{-.55}{0}{agreen}
    \bit{2.35}{-.55}{1}{bred}\bit{2.65}{-.55}{0}{bred}\bit{2.95}{-.55}{0}{bred}\bit{3.25}{-.55}{0}{bred}
  }

  \pathpanel{0}{0}{bred}{\textcolor{agreen}{34}\textcolor{bred}{78}}{$7,8$ chosen}{
    \draw[bred,line width=1.66416pt] (0,0)--(1,0)--(2,0)--(2,1)--(2,2)--(3,2)--(4,2)--(4,3)--(4,4);
    \fill (0,0) circle(1.4pt); \fill (4,4) circle(1.4pt);
    \fill[white] (2,2) circle(2.1pt); \draw (2,2) circle(2.1pt);
    \stepnum{.5}{.18}{1}{agreen}\stepnum{1.5}{.18}{2}{agreen}\stepnum{2.18}{.5}{3}{agreen}\stepnum{2.18}{1.5}{4}{agreen}
    \stepnum{2.5}{2.18}{5}{bred}\stepnum{3.5}{2.18}{6}{bred}\stepnum{3.82}{2.5}{7}{bred}\stepnum{3.82}{3.5}{8}{bred}
    \bit{1.0}{-.55}{0}{agreen}\bit{1.3}{-.55}{0}{agreen}\bit{1.6}{-.55}{1}{agreen}\bit{1.9}{-.55}{1}{agreen}
    \bit{2.35}{-.55}{0}{bred}\bit{2.65}{-.55}{0}{bred}\bit{2.95}{-.55}{1}{bred}\bit{3.25}{-.55}{1}{bred}
  }

  \pathpanel{5.05}{0}{bred}{\textcolor{agreen}{4}\textcolor{bred}{567}}{one $A$-point}{
    \draw[bred,line width=1.66416pt] (0,0)--(1,0)--(2,0)--(3,0)--(3,1)--(3,2)--(3,3)--(3,4)--(4,4);
    \fill (0,0) circle(1.4pt); \fill (4,4) circle(1.4pt);
    \fill[white] (3,1) circle(2.1pt); \draw (3,1) circle(2.1pt);
    \stepnum{.5}{.18}{1}{agreen}\stepnum{1.5}{.18}{2}{agreen}\stepnum{2.5}{.18}{3}{agreen}\stepnum{3.18}{.5}{4}{agreen}
    \stepnum{3.18}{1.5}{5}{bred}\stepnum{3.18}{2.5}{6}{bred}\stepnum{3.18}{3.5}{7}{bred}\stepnum{3.5}{3.82}{8}{bred}
    \bit{1.0}{-.55}{0}{agreen}\bit{1.3}{-.55}{0}{agreen}\bit{1.6}{-.55}{0}{agreen}\bit{1.9}{-.55}{1}{agreen}
    \bit{2.35}{-.55}{1}{bred}\bit{2.65}{-.55}{1}{bred}\bit{2.95}{-.55}{1}{bred}\bit{3.25}{-.55}{0}{bred}
  }

  \pathpanel{10.1}{0}{bred}{\textcolor{agreen}{4}\textcolor{bred}{678}}{both mechanisms}{
    \draw[bred,line width=1.66416pt] (0,0)--(1,0)--(2,0)--(3,0)--(3,1)--(4,1)--(4,2)--(4,3)--(4,4);
    \fill (0,0) circle(1.4pt); \fill (4,4) circle(1.4pt);
    \fill[white] (3,1) circle(2.1pt); \draw (3,1) circle(2.1pt);
    \stepnum{.5}{.18}{1}{agreen}\stepnum{1.5}{.18}{2}{agreen}\stepnum{2.5}{.18}{3}{agreen}\stepnum{3.18}{.5}{4}{agreen}
    \stepnum{3.5}{1.18}{5}{bred}\stepnum{3.82}{1.5}{6}{bred}\stepnum{3.82}{2.5}{7}{bred}\stepnum{3.82}{3.5}{8}{bred}
    \bit{1.0}{-.55}{0}{agreen}\bit{1.3}{-.55}{0}{agreen}\bit{1.6}{-.55}{0}{agreen}\bit{1.9}{-.55}{1}{agreen}
    \bit{2.35}{-.55}{0}{bred}\bit{2.65}{-.55}{1}{bred}\bit{2.95}{-.55}{1}{bred}\bit{3.25}{-.55}{1}{bred}
  }

\end{tikzpicture}
\caption{For $k=3$, examples of $A$-leaning sets (top) and $B$-leaning sets (bottom), encoded as lattice paths. Here $A=1234$ and $B=5678$; step $i$ goes up when $i$ is chosen and right when it is omitted.}
\label{fig:leaning-paths}
\end{figure}
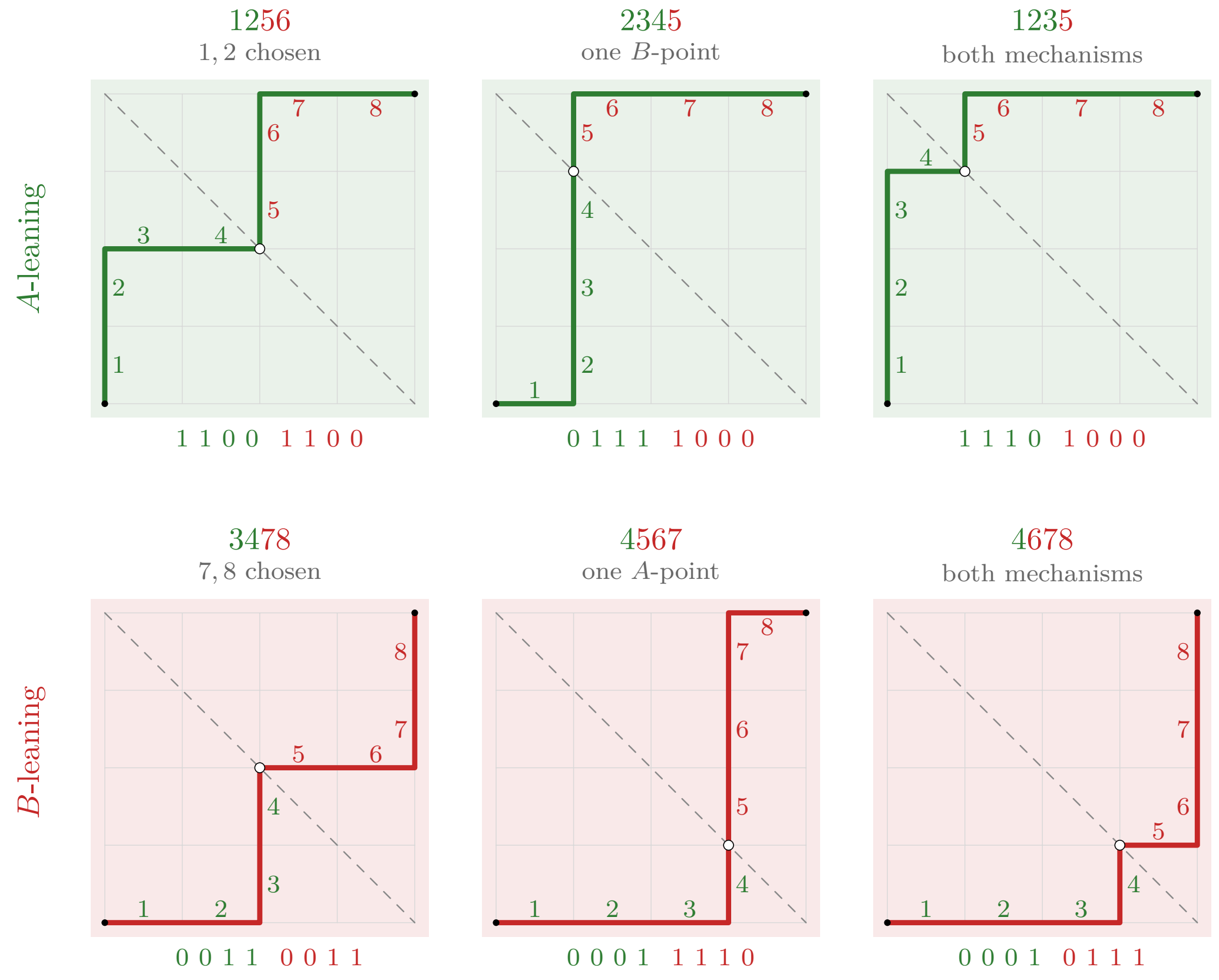

\begin{lemma}\label{lem:local-consistent}
The families $X_{AB}$ and $Y_{AB}$ are consistent.
\end{lemma}

\begin{proof}
$A$ is the first member of every packet containing it, and $B$ is
the last member of every packet containing it. No packet contains both.
So $X_{AB}$ is consistent.

To prove consistency of $Y_{AB}$, add $A$ to the $A$-leaning family. The resulting family is closed under
moving coordinates left. Indeed, moving left cannot introduce $b_k$ or
$b_{k+1}$; a set with one point in $B$ keeps at most one, and a set containing
$a_1,a_2$ keeps these two smallest points. Since packet members increase
coordinatewise in packet order, this extended A-leaning family meets every packet in a beginning segment.

Removing $A$ changes only packets $P(A\cup\{b_j\})$. For $j\le k-1$, every
packet member belonged to the extended A-leaning family; removing its first member $A$ leaves an ending segment. For $j\in\{k,k+1\}$, the packet contains no $A$-leaning members. Thus the $A$-leaning family is consistent. Reflection proves the same
for the $B$-leaning family.

No packet meets both leaning families. If a parent $(k+2)$-set $H$ has a
$B$-leaning facet, it contains at most one of $a_1,a_2$. An $A$-leaning facet
of the same parent must therefore have only one point in $B$, so $H$ has at
most two points in $B$. The reflected argument gives at most two points in
$A$. This contradicts $|H|=k+2\ge5$. The union of the two leaning families is
consequently consistent, and so is its complement $Y_{AB}$.

\end{proof}

\begin{figure}[H]
\centering
\begin{tikzpicture}[scale=.425833,font=\diagramfont,line cap=round,line join=round]
\useasboundingbox (0pt,0pt) rectangle (594.000pt,594.000pt);
\draw[line width=.91980pt] (269.185pt,228.157pt) -- (231.257pt,134.279pt) -- (226.505pt,136.324pt) -- (221.850pt,138.576pt) -- (217.299pt,141.027pt) -- (212.861pt,143.670pt) -- (208.543pt,146.496pt) -- (204.351pt,149.495pt) -- (200.291pt,152.660pt) -- (196.370pt,155.979pt) -- (192.591pt,159.446pt) -- (188.958pt,163.049pt) -- (185.476pt,166.781pt) -- (182.147pt,170.631pt) -- (178.974pt,174.590pt) -- (175.958pt,178.651pt) -- (173.100pt,182.802pt) -- (170.402pt,187.036pt) -- (167.863pt,191.345pt) -- (165.483pt,195.718pt) -- (163.261pt,200.150pt) -- (161.195pt,204.630pt) -- (159.285pt,209.153pt) -- (157.527pt,213.710pt) -- (155.918pt,218.296pt) -- (154.457pt,222.902pt) -- (153.139pt,227.524pt) -- (151.961pt,232.154pt) -- (150.918pt,236.789pt) -- (150.007pt,241.423pt) -- (149.223pt,246.052pt) -- (148.562pt,250.672pt) -- (148.018pt,255.279pt) -- (147.586pt,259.870pt) -- (147.262pt,264.443pt) -- (147.041pt,268.996pt) -- (146.917pt,273.527pt) -- (146.885pt,278.036pt) -- (146.940pt,282.522pt) -- (147.077pt,286.984pt) -- (147.292pt,291.423pt) -- (147.579pt,295.841pt) -- (147.933pt,300.238pt) -- (148.351pt,304.617pt) -- (148.828pt,308.979pt) -- (149.361pt,313.329pt) -- (149.945pt,317.667pt) -- (223.473pt,307.334pt);
\node[font=\diagramfont] at (208.925pt,246.150pt) {$C_{1}$};
\node[anchor=base,font=\diagramfont,text=agreen] at (108.008pt,197.550pt) {1};
\node[anchor=base,font=\diagramfont,text=agreen] at (119.168pt,197.550pt) {2};
\node[anchor=base,font=\diagramfont,text=agreen] at (130.328pt,197.550pt) {3};
\node[anchor=base,font=\diagramfont,text=bred] at (141.488pt,197.550pt) {8};
\draw[line width=.91980pt] (223.473pt,286.666pt) -- (123.208pt,272.575pt) -- (122.603pt,277.712pt) -- (122.226pt,282.870pt) -- (122.073pt,288.037pt) -- (122.143pt,293.202pt) -- (122.431pt,298.354pt) -- (122.933pt,303.484pt) -- (123.643pt,308.582pt) -- (124.557pt,313.638pt) -- (125.670pt,318.644pt) -- (126.974pt,323.592pt) -- (128.465pt,328.473pt) -- (130.135pt,333.281pt) -- (131.977pt,338.009pt) -- (133.985pt,342.651pt) -- (136.152pt,347.201pt) -- (138.470pt,351.655pt) -- (140.931pt,356.008pt) -- (143.529pt,360.256pt) -- (146.256pt,364.396pt) -- (149.103pt,368.425pt) -- (152.065pt,372.341pt) -- (155.132pt,376.143pt) -- (158.299pt,379.828pt) -- (161.558pt,383.397pt) -- (164.901pt,386.849pt) -- (168.322pt,390.185pt) -- (171.815pt,393.405pt) -- (175.373pt,396.511pt) -- (178.990pt,399.505pt) -- (182.659pt,402.387pt) -- (186.377pt,405.162pt) -- (190.138pt,407.831pt) -- (193.936pt,410.398pt) -- (197.769pt,412.866pt) -- (201.631pt,415.239pt) -- (205.519pt,417.521pt) -- (209.431pt,419.716pt) -- (213.364pt,421.829pt) -- (217.316pt,423.863pt) -- (221.285pt,425.823pt) -- (225.271pt,427.715pt) -- (229.272pt,429.542pt) -- (233.289pt,431.310pt) -- (237.321pt,433.024pt) -- (241.371pt,434.687pt) -- (269.185pt,365.843pt);
\node[font=\diagramfont] at (208.925pt,347.850pt) {$C_{2}$};
\node[anchor=base,font=\diagramfont,text=agreen] at (108.008pt,396.450pt) {1};
\node[anchor=base,font=\diagramfont,text=agreen] at (119.168pt,396.450pt) {3};
\node[anchor=base,font=\diagramfont,text=agreen] at (130.328pt,396.450pt) {4};
\node[anchor=base,font=\diagramfont,text=bred] at (141.488pt,396.450pt) {8};
\draw[line width=.91980pt] (251.287pt,355.510pt) -- (188.951pt,435.296pt) -- (193.098pt,438.389pt) -- (197.376pt,441.294pt) -- (201.774pt,444.009pt) -- (206.282pt,446.531pt) -- (210.888pt,448.858pt) -- (215.582pt,450.989pt) -- (220.352pt,452.922pt) -- (225.188pt,454.659pt) -- (230.079pt,456.198pt) -- (235.016pt,457.542pt) -- (239.989pt,458.692pt) -- (244.987pt,459.650pt) -- (250.003pt,460.418pt) -- (255.027pt,461.000pt) -- (260.052pt,461.399pt) -- (265.068pt,461.619pt) -- (270.068pt,461.663pt) -- (275.046pt,461.538pt) -- (279.995pt,461.247pt) -- (284.908pt,460.795pt) -- (289.780pt,460.188pt) -- (294.606pt,459.432pt) -- (299.381pt,458.533pt) -- (304.101pt,457.495pt) -- (308.762pt,456.326pt) -- (313.362pt,455.031pt) -- (317.897pt,453.616pt) -- (322.366pt,452.088pt) -- (326.766pt,450.452pt) -- (331.098pt,448.715pt) -- (335.360pt,446.883pt) -- (339.551pt,444.961pt) -- (343.674pt,442.955pt) -- (347.727pt,440.870pt) -- (351.714pt,438.712pt) -- (355.634pt,436.485pt) -- (359.491pt,434.195pt) -- (363.287pt,431.845pt) -- (367.024pt,429.439pt) -- (370.707pt,426.982pt) -- (374.338pt,424.477pt) -- (377.921pt,421.925pt) -- (381.460pt,419.331pt) -- (384.961pt,416.695pt) -- (388.426pt,414.020pt) -- (342.713pt,355.510pt);
\node[font=\diagramfont] at (297.000pt,398.700pt) {$C_{3}$};
\node[anchor=base,font=\diagramfont,text=agreen] at (280.260pt,495.900pt) {1};
\node[anchor=base,font=\diagramfont,text=agreen] at (291.420pt,495.900pt) {4};
\node[anchor=base,font=\diagramfont,text=bred] at (302.580pt,495.900pt) {7};
\node[anchor=base,font=\diagramfont,text=bred] at (313.740pt,495.900pt) {8};
\draw[line width=.91980pt] (324.815pt,365.843pt) -- (362.743pt,459.721pt) -- (367.495pt,457.676pt) -- (372.150pt,455.424pt) -- (376.701pt,452.973pt) -- (381.139pt,450.330pt) -- (385.457pt,447.504pt) -- (389.649pt,444.505pt) -- (393.709pt,441.340pt) -- (397.630pt,438.021pt) -- (401.409pt,434.554pt) -- (405.042pt,430.951pt) -- (408.524pt,427.219pt) -- (411.853pt,423.369pt) -- (415.026pt,419.410pt) -- (418.042pt,415.349pt) -- (420.900pt,411.198pt) -- (423.598pt,406.964pt) -- (426.137pt,402.655pt) -- (428.517pt,398.282pt) -- (430.739pt,393.850pt) -- (432.805pt,389.370pt) -- (434.715pt,384.847pt) -- (436.473pt,380.290pt) -- (438.082pt,375.704pt) -- (439.543pt,371.098pt) -- (440.861pt,366.476pt) -- (442.039pt,361.846pt) -- (443.082pt,357.211pt) -- (443.993pt,352.577pt) -- (444.777pt,347.948pt) -- (445.438pt,343.328pt) -- (445.982pt,338.721pt) -- (446.414pt,334.130pt) -- (446.738pt,329.557pt) -- (446.959pt,325.004pt) -- (447.083pt,320.473pt) -- (447.115pt,315.964pt) -- (447.060pt,311.478pt) -- (446.923pt,307.016pt) -- (446.708pt,302.577pt) -- (446.421pt,298.159pt) -- (446.067pt,293.762pt) -- (445.649pt,289.383pt) -- (445.172pt,285.021pt) -- (444.639pt,280.671pt) -- (444.055pt,276.333pt) -- (370.527pt,286.666pt);
\node[font=\diagramfont] at (385.075pt,347.850pt) {$C_{4}$};
\node[anchor=base,font=\diagramfont,text=agreen] at (452.512pt,396.450pt) {1};
\node[anchor=base,font=\diagramfont,text=bred] at (463.672pt,396.450pt) {6};
\node[anchor=base,font=\diagramfont,text=bred] at (474.832pt,396.450pt) {7};
\node[anchor=base,font=\diagramfont,text=bred] at (485.992pt,396.450pt) {8};
\draw[line width=.91980pt] (370.527pt,307.334pt) -- (470.792pt,321.425pt) -- (471.397pt,316.288pt) -- (471.774pt,311.130pt) -- (471.927pt,305.963pt) -- (471.857pt,300.798pt) -- (471.569pt,295.646pt) -- (471.067pt,290.516pt) -- (470.357pt,285.418pt) -- (469.443pt,280.362pt) -- (468.330pt,275.356pt) -- (467.026pt,270.408pt) -- (465.535pt,265.527pt) -- (463.865pt,260.719pt) -- (462.023pt,255.991pt) -- (460.015pt,251.349pt) -- (457.848pt,246.799pt) -- (455.530pt,242.345pt) -- (453.069pt,237.992pt) -- (450.471pt,233.744pt) -- (447.744pt,229.604pt) -- (444.897pt,225.575pt) -- (441.935pt,221.659pt) -- (438.868pt,217.857pt) -- (435.701pt,214.172pt) -- (432.442pt,210.603pt) -- (429.099pt,207.151pt) -- (425.678pt,203.815pt) -- (422.185pt,200.595pt) -- (418.627pt,197.489pt) -- (415.010pt,194.495pt) -- (411.341pt,191.613pt) -- (407.623pt,188.838pt) -- (403.862pt,186.169pt) -- (400.064pt,183.602pt) -- (396.231pt,181.134pt) -- (392.369pt,178.761pt) -- (388.481pt,176.479pt) -- (384.569pt,174.284pt) -- (380.636pt,172.171pt) -- (376.684pt,170.137pt) -- (372.715pt,168.177pt) -- (368.729pt,166.285pt) -- (364.728pt,164.458pt) -- (360.711pt,162.690pt) -- (356.679pt,160.976pt) -- (352.629pt,159.313pt) -- (324.815pt,228.157pt);
\node[font=\diagramfont] at (385.075pt,246.150pt) {$C_{5}$};
\node[anchor=base,font=\diagramfont,text=agreen] at (452.512pt,197.550pt) {1};
\node[anchor=base,font=\diagramfont,text=bred] at (463.672pt,197.550pt) {5};
\node[anchor=base,font=\diagramfont,text=bred] at (474.832pt,197.550pt) {6};
\node[anchor=base,font=\diagramfont,text=bred] at (485.992pt,197.550pt) {8};
\draw[line width=.91980pt] (342.713pt,238.490pt) -- (405.049pt,158.704pt) -- (400.902pt,155.611pt) -- (396.624pt,152.706pt) -- (392.226pt,149.991pt) -- (387.718pt,147.469pt) -- (383.112pt,145.142pt) -- (378.418pt,143.011pt) -- (373.648pt,141.078pt) -- (368.812pt,139.341pt) -- (363.921pt,137.802pt) -- (358.984pt,136.458pt) -- (354.011pt,135.308pt) -- (349.013pt,134.350pt) -- (343.997pt,133.582pt) -- (338.973pt,133.000pt) -- (333.948pt,132.601pt) -- (328.932pt,132.381pt) -- (323.932pt,132.337pt) -- (318.954pt,132.462pt) -- (314.005pt,132.753pt) -- (309.092pt,133.205pt) -- (304.220pt,133.812pt) -- (299.394pt,134.568pt) -- (294.619pt,135.467pt) -- (289.899pt,136.505pt) -- (285.238pt,137.674pt) -- (280.638pt,138.969pt) -- (276.103pt,140.384pt) -- (271.634pt,141.912pt) -- (267.234pt,143.548pt) -- (262.902pt,145.285pt) -- (258.640pt,147.117pt) -- (254.449pt,149.039pt) -- (250.326pt,151.045pt) -- (246.273pt,153.130pt) -- (242.286pt,155.288pt) -- (238.366pt,157.515pt) -- (234.509pt,159.805pt) -- (230.713pt,162.155pt) -- (226.976pt,164.561pt) -- (223.293pt,167.018pt) -- (219.662pt,169.523pt) -- (216.079pt,172.075pt) -- (212.540pt,174.669pt) -- (209.039pt,177.305pt) -- (205.574pt,179.980pt) -- (251.287pt,238.490pt);
\node[font=\diagramfont] at (297.000pt,195.300pt) {$C_{6}$};
\node[anchor=base,font=\diagramfont,text=agreen] at (280.260pt,98.100pt) {1};
\node[anchor=base,font=\diagramfont,text=agreen] at (291.420pt,98.100pt) {2};
\node[anchor=base,font=\diagramfont,text=bred] at (302.580pt,98.100pt) {5};
\node[anchor=base,font=\diagramfont,text=bred] at (313.740pt,98.100pt) {8};
\draw[-{Triangle[length=1.91625pt,width=1.91625pt]},line width=.68985pt] (90.215pt,225.798pt) .. controls (49.572pt,273.266pt) and (49.572pt,320.734pt) .. (90.215pt,368.202pt);
\filldraw[fill=white,draw=black,line width=.55358pt] (59.733pt,297.000pt) ellipse [x radius=54.000pt,y radius=21.000pt];
\path[fill=lightgreen] (29.333pt,287.000pt) rectangle (81.053pt,310.500pt);
\node[anchor=base,font=\diagramdetailfont,text=agreen] at (33.533pt,290.200pt) {1};
\node[anchor=base,font=\diagramdetailfont,text=agreen] at (44.933pt,290.200pt) {2};
\node[anchor=base,font=\diagramdetailfont,text=agreen] at (56.333pt,290.200pt) {3};
\node[anchor=base,font=\diagramdetailfont,text=agreen] at (67.733pt,290.200pt) {4};
\node[anchor=base,font=\diagramdetailfont,text=bred] at (85.933pt,290.200pt) {8};
\draw[-{Triangle[length=1.91625pt,width=1.91625pt]},line width=.68985pt] (131.945pt,440.480pt) .. controls (152.732pt,499.412pt) and (193.840pt,523.146pt) .. (255.270pt,511.682pt);
\filldraw[fill=white,draw=black,line width=.55358pt] (178.366pt,502.480pt) ellipse [x radius=54.000pt,y radius=21.000pt];
\path[fill=lightred] (166.166pt,492.480pt) rectangle (217.886pt,515.980pt);
\node[anchor=base,font=\diagramdetailfont,text=agreen] at (152.166pt,495.680pt) {1};
\node[anchor=base,font=\diagramdetailfont,text=agreen] at (170.366pt,495.680pt) {3};
\node[anchor=base,font=\diagramdetailfont,text=agreen] at (181.766pt,495.680pt) {4};
\node[anchor=base,font=\diagramdetailfont,text=bred] at (193.166pt,495.680pt) {7};
\node[anchor=base,font=\diagramdetailfont,text=bred] at (204.566pt,495.680pt) {8};
\draw[-{Triangle[length=1.91625pt,width=1.91625pt]},line width=.68985pt] (338.730pt,511.682pt) .. controls (400.160pt,523.146pt) and (441.268pt,499.412pt) .. (462.055pt,440.480pt);
\filldraw[fill=white,draw=black,line width=.55358pt] (415.634pt,502.480pt) ellipse [x radius=54.000pt,y radius=21.000pt];
\path[fill=lightred] (403.434pt,492.480pt) rectangle (455.154pt,515.980pt);
\node[anchor=base,font=\diagramdetailfont,text=agreen] at (389.434pt,495.680pt) {1};
\node[anchor=base,font=\diagramdetailfont,text=agreen] at (407.634pt,495.680pt) {4};
\node[anchor=base,font=\diagramdetailfont,text=bred] at (419.034pt,495.680pt) {6};
\node[anchor=base,font=\diagramdetailfont,text=bred] at (430.434pt,495.680pt) {7};
\node[anchor=base,font=\diagramdetailfont,text=bred] at (441.834pt,495.680pt) {8};
\draw[-{Triangle[length=1.91625pt,width=1.91625pt]},line width=.68985pt] (503.785pt,368.202pt) .. controls (544.428pt,320.734pt) and (544.428pt,273.266pt) .. (503.785pt,225.798pt);
\filldraw[fill=white,draw=black,line width=.55358pt] (534.267pt,297.000pt) ellipse [x radius=54.000pt,y radius=21.000pt];
\path[fill=lightred] (522.067pt,287.000pt) rectangle (573.787pt,310.500pt);
\node[anchor=base,font=\diagramdetailfont,text=agreen] at (508.067pt,290.200pt) {1};
\node[anchor=base,font=\diagramdetailfont,text=bred] at (526.267pt,290.200pt) {5};
\node[anchor=base,font=\diagramdetailfont,text=bred] at (537.667pt,290.200pt) {6};
\node[anchor=base,font=\diagramdetailfont,text=bred] at (549.067pt,290.200pt) {7};
\node[anchor=base,font=\diagramdetailfont,text=bred] at (560.467pt,290.200pt) {8};
\draw[-{Triangle[length=1.91625pt,width=1.91625pt]},line width=.68985pt] (462.055pt,153.520pt) .. controls (441.268pt,94.588pt) and (400.160pt,70.854pt) .. (338.730pt,82.318pt);
\filldraw[fill=white,draw=black,line width=.55358pt] (415.634pt,91.520pt) ellipse [x radius=54.000pt,y radius=21.000pt];
\path[fill=lightgreen] (385.234pt,81.520pt) rectangle (436.954pt,105.020pt);
\node[anchor=base,font=\diagramdetailfont,text=agreen] at (389.434pt,84.720pt) {1};
\node[anchor=base,font=\diagramdetailfont,text=agreen] at (400.834pt,84.720pt) {2};
\node[anchor=base,font=\diagramdetailfont,text=bred] at (412.234pt,84.720pt) {5};
\node[anchor=base,font=\diagramdetailfont,text=bred] at (423.634pt,84.720pt) {6};
\node[anchor=base,font=\diagramdetailfont,text=bred] at (441.834pt,84.720pt) {8};
\draw[-{Triangle[length=1.91625pt,width=1.91625pt]},line width=.68985pt] (255.270pt,82.318pt) .. controls (193.840pt,70.854pt) and (152.732pt,94.588pt) .. (131.945pt,153.520pt);
\filldraw[fill=white,draw=black,line width=.55358pt] (178.366pt,91.520pt) ellipse [x radius=54.000pt,y radius=21.000pt];
\path[fill=lightgreen] (147.966pt,81.520pt) rectangle (199.686pt,105.020pt);
\node[anchor=base,font=\diagramdetailfont,text=agreen] at (152.166pt,84.720pt) {1};
\node[anchor=base,font=\diagramdetailfont,text=agreen] at (163.566pt,84.720pt) {2};
\node[anchor=base,font=\diagramdetailfont,text=agreen] at (174.966pt,84.720pt) {3};
\node[anchor=base,font=\diagramdetailfont,text=bred] at (186.366pt,84.720pt) {5};
\node[anchor=base,font=\diagramdetailfont,text=bred] at (204.566pt,84.720pt) {8};
\fill[white] (239.625pt,396.376pt) circle[radius=14.850pt];
\node[font=\diagramfont,text=agreen] at (239.625pt,396.376pt) {4};
\fill[white] (354.375pt,396.376pt) circle[radius=14.850pt];
\node[font=\diagramfont,text=bred] at (354.375pt,396.376pt) {7};
\fill[white] (411.750pt,297.000pt) circle[radius=14.850pt];
\node[font=\diagramfont,text=bred] at (411.750pt,297.000pt) {6};
\fill[white] (354.375pt,197.624pt) circle[radius=14.850pt];
\node[font=\diagramfont,text=bred] at (354.375pt,197.624pt) {5};
\fill[white] (239.625pt,197.624pt) circle[radius=14.850pt];
\node[font=\diagramfont,text=agreen] at (239.625pt,197.624pt) {2};
\fill[white] (182.250pt,297.000pt) circle[radius=14.850pt];
\node[font=\diagramfont,text=agreen] at (182.250pt,297.000pt) {3};
\node[font=\diagramfont] at (297.000pt,297.000pt)
  {\textcolor{agreen}{1}\quad\textcolor{bred}{8}};
\end{tikzpicture}
\caption{The local obstruction for $k=3$. The six sets $C_1,\ldots,C_6$ form the forced cycle. The outer bubbles show the packet parents; the pale green or red rectangles mark the packet members which are $A$, $A$-leaning, $B$ or $B$-leaning -- these fix the packet directions.}
\label{fig:ziegler-cycle}
\end{figure}
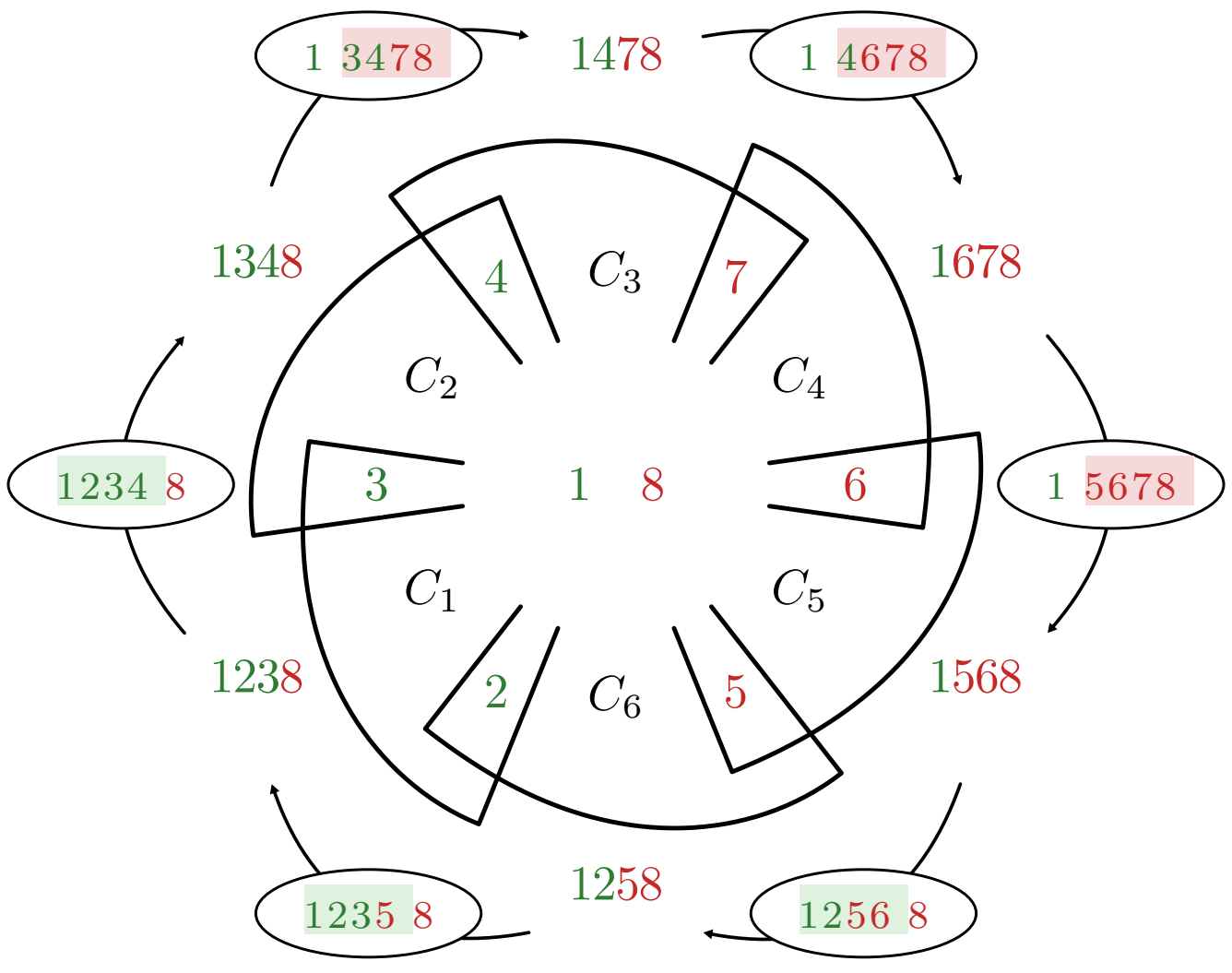

\begin{lemma}\label{lem:local-gap}
The pair $(X_{AB},Y_{AB})$ is a gap.
\end{lemma}

\begin{proof}
Only mixed sets are leaning, so $X_{AB}\subseteq Y_{AB}$. Read
\[
a_2,a_3,\ldots,a_{k+1},b_k,b_{k-1},\ldots,b_1
\]
cyclically. Starting at each of its $2k$ positions, take $k-1$ consecutive
points and add $a_1,b_{k+1}$. Call the resulting sets
$C_1,\ldots,C_{2k}$ in cyclic order. They contain $a_1,b_{k+1}$, hence belong to
$Y_{AB}\setminus X_{AB}$. In particular,
\[
C_1=A-a_{k+1}+b_{k+1},\quad C_2=A-a_2+b_{k+1},
\]
\[
C_{k+1}=B-b_1+a_1,\quad C_{k+2}=B-b_k+a_1.
\]
Assume a flip path from $X_{AB}$ to $Y_{AB}$ exists, and write $F\prec G$
when $F$ is added before $G$. Every family on the path lies between
$X_{AB}$ and $Y_{AB}$ under inclusion.

Every mixed facet of $P(A\cup\{b_{k+1}\})$ contains $b_{k+1}$ and at least one
of $a_1,a_2$. Thus this packet is contained in $Y_{AB}$. Initially only its
first member $A$ is present. Its other members must be added from first to
last, giving $C_1\prec C_2$.

For $2\le i\le k$, the parent $C_i\cup C_{i+1}$ contains $a_1$ but not
$a_2$, and contains both $b_k,b_{k+1}$. Its only facet outside $Y_{AB}$ is
obtained by deleting $a_1$, the last packet member. Every nonempty
intersection with this packet along the proposed path must therefore be a
beginning segment. Passing from $C_i$ to $C_{i+1}$ replaces an $A$-point by a
larger $B$-point, so $C_i$ precedes $C_{i+1}$ in the packet. Hence
\[
C_1\prec C_2\prec\cdots\prec C_{k+1}.
\]
Every mixed facet of $P(\{a_1\}\cup B)$ contains $a_1$ and at least one of
$b_k,b_{k+1}$, so this packet is also contained in $Y_{AB}$. Initially only
its last member $B$ is present, giving $C_{k+1}\prec C_{k+2}$.

For $k+2\le i\le2k$, with $C_{2k+1}=C_1$, the parent
$C_i\cup C_{i+1}$ contains $a_1,a_2,b_{k+1}$ but not $b_k$. Its only facet
outside $Y_{AB}$ is obtained by deleting $b_{k+1}$, the first packet member.
Every nonempty intersection along the path must therefore be an ending
segment. Passing from $C_i$ to $C_{i+1}$ replaces a $B$-point by a smaller
$A$-point. Thus $C_i$ follows $C_{i+1}$ in packet order and must be added
first. We obtain
\[
C_1\prec C_2\prec\cdots\prec C_{2k}\prec C_1,
\]
a contradiction.
\end{proof}

\subsection{Globalization}
Let $L<R$ be consecutive blocks, each of size at least $k+1$, and take
$A\in\binom{L}{k+1}$ and $B\in\binom{R}{k+1}$. Call a set left-pure if it
lies in $L$, right-pure if it lies in $R$, and mixed if it meets both. For sets $F$ and $G$ of equal cardinality, write $F \preceq G$ when $f_i \leq g_i$ for every $i$. We now
construct consistent families $x_{AB},y_{AB}$ whose restriction to $A\cup B$
is the local gap of the previous section. We will then build a standard example of many such gaps.

\begin{definition}\label{def:x}
\[
x_{AB}=\left\{F\in\binom{L}{k+1}:F\preceq A\right\}
\cup
\left\{F\in\binom{R}{k+1}:B\preceq F\right\}.
\]
\end{definition}

\begin{lemma}\label{x-cons}
    $x_{AB}$ is consistent and $x_{AB}|_{A \cup B} = X_{AB}$.
\end{lemma}

\begin{proof}
On a packet whose parent lies in $L$, the intersection with $\left\{F\in\binom{L}{k+1}:F\preceq A\right\}$ is a
beginning segment, because $\left\{F\in\binom{L}{k+1}:F\preceq A\right\}$ is a down-set in $\preceq$-order. Similarly, on a packet whose parent lies in $R$, the intersection
with $\left\{F\in\binom{R}{k+1}:B\preceq F\right\}$ is an ending segment. A mixed packet parent has at most one pure facet,
first if left-pure and last if right-pure, and $x_{AB}$ consists only of pure sets. Thus $x_{AB}$ is consistent. Its restriction to $A \cup B$ is $\{A, B\} = X_{AB}$.
\end{proof}

To define $y_{AB}$, we will use an auxiliary construction of {\em trails}, motivated by the construction of operadic insertions in \cite[Section~5.8]{BruhatOperads}.

\begin{definition}\label{def:trail}
For $P\subseteq L$, match its points from right to left to the largest unused
targets in $A$ not greater than them. If the matching succeeds, denote the set of matched targets
by $\tau_A(P)$, the left trail; otherwise say that the left trail is undefined. For $Q\subseteq R$, match its points
from left to right to the smallest unused targets in $B$ not smaller
than them. If the matching succeeds, denote the set of matched targets
by $\tau_B(Q)$, the right trail; otherwise say that the right trail is undefined. For $S\subseteq L\cup R$, define
\[
\tau_{AB}(S)=\tau_A(S\cap L)\cup\tau_B(S\cap R)
\]
if both sides exist; otherwise say that the total trail is undefined.
\end{definition}

In other words, the left trail $\tau_A(P)$ is the greatest $|P|$-subset of $A$ that is
coordinatewise at most $P$; the right trail $\tau_B(Q)$ is the least $|Q|$-subset of
$B$ that is coordinatewise at least $Q$. 

\begin{definition}
    $y_{AB}$ consists of sets in $\binom{L\cup R}{k+1}$ whose total trail is undefined or belongs to $Y_{AB}$.
\end{definition}

\begin{lemma}
    $y_{AB}$ is consistent and $y_{AB}|_{A \cup B} = Y_{AB}$.
\end{lemma}

\begin{proof}
Let $H \subset L \cup R$ be a packet parent, with $P = H \cap L$ and $Q = H \cap R$. If both $\tau_A(P)$ and $\tau_B(Q)$ are undefined, then removing one point from $H$ can only repair one of the two trails. Thus every facet of $H$ has undefined total trail and the whole packet is in $y_{AB}$. Now assume that only $\tau_{A}(P)$ is undefined. Denote the points of $A$ by $a_1< \ldots < a_{k+1}$, put $a_{k+2}=+\infty$, and denote the points of $P$ by $p_1<\ldots < p_s$. The left trail is defined exactly when $a_i\leq p_i$ for $1\leq i\leq s$, so let $d$ be the smallest index where $a_d>p_d$. Deleting points of $Q$ or points $p_i$ for $i>d$ does not repair the left trail. The deletions $p_i$ with $i \leq d$ either all fail to repair it, in which case the whole packet lies in $y_{AB}$, or all give one and the same total trail $T$. In the latter case, if this $T$ is in $Y_{AB}$, then the whole packet of $H$ is in $y_{AB}$, and otherwise the intersection is a beginning segment. The case when only $\tau_{B}(Q)$ is undefined is dual. Finally, assume that both $\tau_{A}(P)$ and $\tau_{B}(Q)$ are defined, and so is $\tau_{AB}(H)$. The trails of the facets of $H$ are facets of $\tau_{AB}(H)$ and occur in its packet order, with repetitions and omissions. The intersection of $Y_{AB}$ with the packet of $\tau_{AB}(H)$ is a beginning or an ending segment by consistency of $Y_{AB}$, so the same follows for $y_{AB}$. Finally, trails fix subsets of $A \cup B$, so $y_{AB}$ restricts to $Y_{AB}$.
\end{proof}

\begin{lemma}\label{lem:global-gap}
The pair $(x_{AB},y_{AB})$ is a gap.
\end{lemma}

\begin{proof} We have already shown that $x_{AB}$ and $y_{AB}$ are consistent. Every pure $(k+1)$-set belongs to $y_{AB}$: its trail is either undefined or equals $A$ or $B$. We have only pure sets in $x_{AB}$, so $x_{AB}\subseteq y_{AB}$. A flip path from $x_{AB}$ to $y_{AB}$ would restrict, after removing
repetitions, to a flip path from $X_{AB}$ to $Y_{AB}$, contradicting that $(X_{AB}, Y_{AB})$ is a gap.
\end{proof}

\subsection{Cross comparisons for different gaps}\label{sec:cross}
We now derive a criterion for cross comparisons between different gaps introduced in the previous subsection. 

\begin{lemma}\label{lem:cross}
For $A,C\in\binom{L}{k+1}$ and $B,D\in\binom{R}{k+1}$,
\begin{equation}\label{eq:cross}
x_{CD}<_{\mathrm{flip}}y_{AB}
\quad\Longleftrightarrow\quad
A\npreceq C\text{ or }D\npreceq B.
\end{equation}
\end{lemma}

\begin{proof}
    First, suppose that $A \preceq C$ and $D \preceq B$. Then $x_{CD}|_{A \cup B} = \{ A, B\} = X_{AB}$, so a flip path from $x_{CD}$ to $y_{AB}$ would restrict, after removing repetitions, to a flip path from $X_{AB}$ to $Y_{AB}$, contradicting that $(X_{AB},Y_{AB})$ is a gap. 

    Now assume $A \npreceq C$ (the case $D \npreceq B$ is symmetric). Introduce the auxiliary family 

    \[
z_{A}=\left\{F\in\binom{L}{k+1}:A \npreceq F\right\}
\cup
\binom{R}{k+1}.
\]

It is consistent by the same argument as we used for $x_{AB}$ in Lemma \ref{x-cons}. We have $x_{CD} \subseteq z_{A}$ because the left-pure part of $x_{CD}$ consists of sets $F \preceq C$, and $A \preceq F$ would imply $A \preceq C$. This inclusion can be realized via a flip path: we can add the missing left-pure sets in the order that leaves the left-pure part a coordinatewise down-set, hence consistent, and then the missing right-pure sets in the order that leaves the right-pure part a coordinatewise up-set, hence also consistent. So $x_{CD}\leq_{\mathrm{flip}} z_{A}$. Since $z_A$ consists only of pure sets, $z_{A} \subseteq y_{AB}$, and it remains to show that the missing sets from $y_{AB}$ can be added while preserving consistency. We do so in the following 3 steps. As {\bf Step 1}, we choose directions for all packets. As {\bf Step 2}, we verify that for every packet, both $z_{A}$ and $y_{AB}$ are beginning segments of the packet in the chosen direction. Finally, as {\bf Step 3}, we demonstrate that the family $W$ of packet parents whose packets are directed in the usual packet order is itself a consistent family of $(k+2)$-sets. Then it would correspond to some admissible ordering of $(k+1)$-sets, and missing elements of $y_{AB}$ will be added to $z_{A}$ in that order.

{\bf Step 1: } We define packet orientations as follows. Consider a packet parent $H$. 
\begin{itemize}
    \item If the total trail $\tau_{AB}(H)$ is defined, orient the packet of $H$ positively if the packet of $\tau_{AB}(H)$ meets $Y_{AB}$ in a nonempty proper beginning segment, and negatively otherwise.
    \item If the left trail $\tau_A(H \cap L)$ is undefined, orient the packet of $H$ positively if $|H \cap L| \geq 2$ and negatively when $|H \cap L| = 1$.
    \item If the left trail $\tau_A(H \cap L)$ is defined but the right trail $\tau_B(H \cap R)$ is undefined, orient $H$ negatively.
\end{itemize}

{\bf Step 2: } We now check that for the chosen orientations of packets, intersections with $z_A$ and $y_{AB}$ are beginning segments. First consider $y_{AB}$. In the first case, where $\tau_{AB}(H)$ is defined, the intersection of its packet with $Y_{AB}$ is a beginning segment in the chosen orientation by the rule, and the same follows for the intersection of the packet of $H$ with $y_{AB}$, as trails of facets follow packet order with possible repetitions and omissions.

Suppose that the left trail is undefined. If $|H\cap L|\geq 2$, then in positive packet order the facets with undefined total trail come first, while any remaining facets have one and the same defined total trail. Thus the intersection with $y_{AB}$ is either the whole packet or a beginning segment. If $|H\cap L|=1$, then every facet belongs to $y_{AB}$. Indeed, deleting a point of $H\cap R$ leaves the left trail undefined, while deleting the unique point of $H\cap L$ gives a right-pure set whose trail is either undefined or equal to $B$. Thus the negative orientation also makes the intersection a beginning segment.

Finally, suppose that the left trail is defined but the right trail is undefined. In positive packet order, the deletions which repair the right trail come first and give one and the same total trail, while all later facets have undefined total trail. Thus the intersection with $y_{AB}$ is either the whole packet or an ending segment, and becomes a beginning segment in the negative orientation.

Now consider $z_A$. If $H\subseteq L$, then its packet is positively oriented, its facets increase in $\preceq$-order, and the facets $F$ satisfying $A\npreceq F$ form a beginning segment. If $H\subseteq R$, then the whole packet belongs to $z_A$.

Suppose that $H$ is mixed. If it has at least two points in both $L$ and $R$, then it has no pure facets and its packet does not intersect $z_A$. If $|H\cap R|=1$, then the only possible member of $z_A$ is the first facet $H\cap L$. If this facet belongs to $z_A$, then $A\npreceq H\cap L$, so its left trail is undefined. Hence the packet of $H$ is positively oriented.

If $|H\cap L|=1$, then the only member of $z_A$ is the last facet $H\cap R$. If the total trail is undefined, the packet is negatively oriented by the second or third rule. If the total trail is defined, then its trail has the form $\{a_i\}\cup B$. Since $B$ is the last member of this packet and belongs to $Y_{AB}$, the intersection with $Y_{AB}$ is either the whole packet or an ending segment. Thus the packet is again negatively oriented, and its intersection with $z_A$ is a beginning segment in the chosen orientation.

{\bf Step 3:} We are left to show that the family $W$ of packet parents whose packets are positively oriented is consistent. Let $K$ be some $(k+3)$-set. First assume that the total trail $\tau_{AB}(K)$ is defined; then total trails of all its facets are also defined and follow the packet order with possible repetitions and omissions. It thus suffices to show that for any $(k+3)$-subset $T$ of $A \cup B$, its positively oriented facets $J$ form a beginning or an ending segment. By our rule, positivity occurs in two cases: 

\begin{enumerate}
    \item $J$ has between $2$ and $k$ points in $A$, contains at most one
    of $a_1,a_2$, and either has exactly two points in $A$ or contains
    $b_k,b_{k+1}$. In this case the packet of $J$ contains a nonempty proper
    ending segment of $B$-leaning sets.
    \item $J=A\cup\{b_j\}$ for some $j\leq k-1$. In this case all packet
    members except the first one, $A$, are $A$-leaning.
\end{enumerate}

Let $c$ be the number of $A$-points in $T$, $2\leq c\leq k+1$. If $c=2$, the facets obtained by deleting points of $T\cap B$ are either
all positive or all negative. The facets obtained by deleting either of
the two $A$-points have only one $A$-point and are negative. Thus
positivity changes at most once.

Suppose that $3\leq c\leq k$. First come the facets obtained by deleting
points of $T\cap B$. Among these, positivity can only change from negative
to positive: deleting $b_{k+1}$ or $b_k$ destroys the condition that both
remain, whereas deleting any smaller $B$-point preserves it. Next come the facets obtained by deleting points of $T\cap A$. Among
these, positivity is either constant or changes from negative to positive
at the final deletions of $a_2$ and $a_1$. Moreover, if the first block
contains a positive facet, then $T$ already contains $b_k,b_{k+1}$ and at
most one of $a_1,a_2$, so every facet in the second block is positive.
Hence the two blocks together have at most one change.

If $c=k+1$, then $T\cap A=A$ and $T\cap B$ has two
points. The first two facets, obtained by deleting these two $B$-points,
have signs \(++\), \(+-\), or \(--\), according to whether the remaining $B$-point lies among
$b_1,\ldots,b_{k-1}$. The later facets are all negative unless \(T\cap B=\{b_k,b_{k+1}\}\). In that exceptional case, precisely the final two facets, obtained by
deleting $a_2$ and $a_1$, are positive, while the first two facets are
negative. Again positivity changes at most once.

Now assume that the left trail of $K\cap L$ is undefined, and put $s=|K\cap L|$. Deleting a point of $K\cap R$ does not repair it. Among the deletions of points of $K\cap L$, the repairing deletions, if there are any, form a final block in packet order and all give the same left trail.

If $s\geq 3$, every facet before this final block still has undefined left trail and at least two left points, so it is positively oriented. The facets in the final block all have the same orientation: if the right trail is defined, they have one and the same total trail, and if it is undefined, they are all negatively oriented by the third rule. Thus positivity changes at most once.

If $s=2$, the facets obtained by deleting points of $K\cap R$ are positively oriented. Every facet obtained by deleting a point of $K\cap L$ is negatively oriented. Indeed, its left part has one point. If its left trail is undefined, this follows from the second rule, and if its left trail is defined while the right trail is undefined, it follows from the third rule. If its total trail is defined, it has one point in $A$ and all of $B$. Since $B$ is the last member of its packet and belongs to $Y_{AB}$, this packet is negatively oriented. Thus the positive facets form a beginning segment. 

If $s=1$, every facet is negatively oriented: deleting a point of $K\cap R$ leaves an undefined left trail with one left point, while deleting the unique point of $K\cap L$ leaves a right-pure $(k+2)$-set whose right trail is undefined.

Finally, suppose that the left trail of $K\cap L$ is defined but the right trail of $K\cap R$ is undefined. Among the deletions of points of $K\cap R$, the repairing deletions form an initial block in packet order and all give the same right trail. Since the left part is unchanged, all facets in this initial block have the same total trail and hence the same orientation. Every remaining facet is negatively oriented. After any remaining deletion from $K\cap R$, the right trail is still undefined, while deleting a point of $K\cap L$ leaves the right part unchanged and the left trail defined. Thus positivity changes at most once.

\end{proof}

\begin{remark}
    The architecture of the proof above bears nontrivial similarity to the architecture of the proof of \cite[Theorem~11.1]{BruhatOperads}, written by the same author.
\end{remark}

\subsection{The standard example}

Order the index pairs by
\[
(A,B)\unlhd(C,D)\quad\Longleftrightarrow\quad
A\preceq C\text{ and }D\preceq B.
\]

\begin{lemma}\label{lem:standard-example}
Any antichain $\mathcal I$ of index pairs gives a
standard example $S_{|\mathcal I|}$ consisting of its global gaps \cite[Chapter~1]{Trotter}. In
particular, $\dim B(n,k)\ge|\mathcal I|$ for every nonempty such antichain,
where $n=|L|+|R|$.
\end{lemma}

\begin{proof}
For such an antichain $\mathcal I$, every
matching pair $(x_{AB},y_{AB})$ is incomparable because it is a gap. For distinct
indices $(A,B),(C,D)$, Lemma~\ref{lem:cross} gives \( x_{CD}<y_{AB} \) and \( x_{AB}<y_{CD} \). 
These relations force all selected elements to be distinct and rule out any
further comparisons. For example, $x_i\le x_j$ would give
$x_i\le x_j<y_i$, a contradiction; the argument for the $y_i$ is dual.
Thus they induce $S_{|\mathcal I|}$. Its dimension gives the lower bound. 
\end{proof}

Identify $L$ and $R$ increasingly with $[\ell]$ and $[m]$, respectively,
and let all sums below refer to these positions. Define
\[
\Delta(A,B)=\sum_{b\in B}b-\sum_{a\in A}a.
\]

\begin{lemma}\label{lem:delta}
Every level set of $\Delta$ is an antichain of index pairs.
\end{lemma}

\begin{proof}
If $(A,B)\unlhd(C,D)$ and the pairs are distinct, then $A\preceq C$ and
$D\preceq B$, with at least one strict comparison. The sum of the left
coordinates increases and the sum of the right coordinates decreases, with
at least one strict change. Hence $\Delta(A,B)>\Delta(C,D)$.
\end{proof}

\begin{proof}[Proof of Theorem~\ref{thm:flip-main}]
Split $[n]$ into consecutive blocks of sizes
$\ell=\lfloor n/2\rfloor$ and $m=\lceil n/2\rceil$. There are
$\binom{\ell}{k+1}\binom{m}{k+1}$ index pairs. The smallest and largest
sums of $(k+1)$-subsets of $[s]$ differ by $(k+1)(s-k-1)$. Consequently,
$\Delta$ has at most
\[
(k+1)(\ell-k-1)+(k+1)(m-k-1)+1
=(k+1)(n-2k-2)+1
\]
values. Its largest level has at least the average number of index pairs.

Finally, $\ell,m=n/2+O(1)$, so
\[
\frac{\binom{\ell}{k+1}\binom{m}{k+1}}
{(k+1)(n-2k-2)+1}
=
\frac{n^{2k+1}}{2^{2k+2}(k+1)((k+1)!)^2}+O_k(n^{2k}).
\]
This finishes the proof.
\end{proof}

\begin{question}
For every fixed $k\geq 3$, does $\dim B(n,k)$ admit a polynomial upper bound in $n$?
\end{question}

{\bf AI disclosure.} ChatGPT Astra was used at all stages of the project, except for asking the questions. Astra herself did not notice that the dimension of the flip order is asymptotically much larger than that of the inclusion order. The human herself would have not noticed the beautiful connection to Tur\'an systems. The human takes full responsibility for the results and finds them interesting.

\end{document}